\documentclass{ifacconf}

\usepackage{graphicx}      
\usepackage{amsfonts}
\usepackage{amsmath}
\usepackage{amssymb}
\usepackage{makeidx}
\usepackage{epsfig}
\usepackage{graphicx}
\usepackage[round, sort, numbers]{natbib}
\begin{document}
\baselineskip 10.5pt
%
%
%
\begin{frontmatter}
\title{\LARGE Safe Human-Inspired Mesoscopic Hybrid Automaton for Longitudinal Vehicle Control}
%
%
%
\author{A. Iovine,}
\author{F. Valentini,}
\author{E. De Santis,}
\author{M. D. Di Benedetto,}
\author{M. Pratesi}
\address{Department of Information Engineering, Computer Science and Mathematics,
 Center of Excellence DEWS, University of L'Aquila.\\(e-mail: \{alessio.iovine\},\{francesco.valentini1\}@graduate.univaq.it, \\
 \{elena.desantis\},\{mariadomenica.dibenedetto\},\{marco.pratesi\}@univaq.it)}
%
%
%
\begin{abstract}                
In this paper a mesoscopic hybrid automaton is introduced in order to obtain a human-inspired based adaptive cruise control. 
The proposed control law fits the design target of replacing and imitating a human driver behaviour. A microscopic hybrid automaton model for longitudinal vehicle control based on human psycho-physical behavior is first presented. 
Then a rule for changing time headway on the basis of macroscopic quantities is used to describe the 
interaction among all next vehicles and their impact on driver performance. Finally, results of the ultimate mesoscopic control model are presented.
\end{abstract}%
\end{frontmatter}%
\section{Introduction}\label{Introduction}%
Traffic control is one of the most studied problems in engineering worldwide.
This is due to its high impact in human life: progress in the knowledge and control of traffic systems would raise life quality (see \cite{StateOFtheART}).
The main goal of traffic control is to improve the traffic management depending on a variety of different goals: congestion, emissions and travel time reduction, safety increments etc.... To this purpose, in the past years a growing development of driver supporting systems took place (e.g. Adaptive Cruise Control (ACC) systems, Advanced Driver-Assistance Systems (ADASs)), which should be able to provide full or partial driver assistance. Once introduced, those tools would need to fit the normal traffic dynamics; therefore, they have to resemble the driver behaviour so that the best option is to let them mimic the human behaviour.

This paper target is to develop an ACC model able to imitate the human way of driving related to comfort while ensuring a proper safety level. 
Over the years, a multitude of them has been generated (see 
\cite{Panwai2005-TITS}). Those models can be classified on the basis of the description level (macroscopic, as in 
\cite{Messmer2000}, microscopic as in \cite{Gazis1961}, \cite{Gipps1981} or mesoscopic, which is a microscopic model that takes into account macroscopic parameters), or the adopted control strategy (centralized, as in \cite{Daganzo}, or decentralized as in \cite{Falconi2012}, \cite{DeSantis2006}).
The focus of this paper is on a decentralized mesoscopic control approach.
We consider $N$ vehicle on a single lane road, sorted by location, indexed by $n\in\left\{1,...,N\right\}$, where $n=1$ denotes the first vehicle on the lane. A hybrid model for each pair $\left(n,n+1\right)$, $n=1,...,N-1$, (the "leader" and the "follower") is developed, based on the classical psycho-physical and stimulus-response car-following models. Our modeling of the traffic flow is therefore microscopic. Since the behaviour of the pair $\left(n,n+1\right)$ in general depends on the behaviour of the pair $\left( n-1,n\right)$, such hybrid systems are interconnected.

We first analyze the properties of the overall hybrid system when applying state feedback control laws, using only local information for each pair of vehicles. Such control laws simulate the human control action 
, where the objective is to minimize the traveling time, while maintaining safety.
In the second part of our work we define a "mesoscopic model", where the control action depends not only on individual information but also on some information about the traffic flow, which is a macroscopic quantity. Such information can be provided by a centralized traffic supervisor, or it can be gathered, elaborated and transmitted by the vehicles themselves, which are supposed to be interconnected not only by the dynamics, but also by a communication network. Thanks to the fact that connected vehicles are nowadays a reality (see \cite{Uhlemann}), we are able to consider this second framework.
The benefits associated with such an information spreading are evaluated in terms of shaving the acceleration peaks reduction and throughput of the highway system. Some simulation results are offered, as well as a discussion on communication questions regarding the proposed model feasibility.

This paper is organized in \ref{Conclusions} sections. 
In Section \ref{Microscopic Hybrid Model} the model of the microscopic hybrid automaton for a single vehicle will be described.
Then in Section \ref{Mesoscopic model} a variance-driven time headway
mechanism will be introduced into the hybrid automaton in order to make it
mesoscopic.
Then Section \ref{vehicular_network} will state available technological solutions for the utilized communication framework, while Section \ref{Simulation results} will provide simulation results about the system behaviour. Conclusions will be offered in Section \ref{Conclusions}.
\section{Microscopic Hybrid Model}\label{Microscopic Hybrid Model}
We have embedded a number of models into a unique model. Since a finite number of control actions have been envisaged 
we think that the hybrid systems framework is the most appropriate one for exactly describing and analyzing the model property. 
The closed loop dynamics of each vehicle is autonomous and affected by a disturbance, which represents the control action of the ahead vehicle. We assume that all vehicles are identical. For a vehicle labelled with $n$, $n+1$ denotes its follower. The
hybrid automaton associated with vehicle $n+1$, with $n=1,...,N-1$, is described
by the tuple%
\begin{equation}
\mathcal{H}_{n+1}=(Q,X,f,Init,Dom,\mathcal{E}%
)\label{HybridAutomaton}%
\end{equation}
where $Q=\left\{  q_{1},..., q_{6}\right\}  $ is the set of discrete
states; $X=\mathbb{R}^{3}$ \ is the continuous state space; $f=\left\{
f_{i},q_{i}\in Q\right\}  $, and $f_{i}:X\times\mathbb{R}\rightarrow
\mathbb{R}^{3}$ \ is a vector field that associates to the discrete state
$q_{i}\in Q $ the continuous time-invariant dynamics
\begin{equation}
\dot{x}(t)=f_{i}(x(t),d(t))\label{diffeq}%
\end{equation}
where $d:\mathbb{R}\rightarrow\mathbb{R}$ is a disturbance; $Init\subseteq
Q\times X$ \ is the set of initial discrete and continuous conditions;
$Dom(\cdot):Q\rightarrow2^{X}$,
 and \ $\mathcal{E}\subseteq Q\times Q$ \ is the set of edges. The automaton hybrid state is the pair $(x,q_{i})\in X\times Q$. Let us
define the function $\boldsymbol{I}:X\rightarrow \left\{ 1,2,...6\right\} $,
where $\boldsymbol{I}(x)=i:x\in Dom(q_{i})$, $Dom(q_{i})\cap Dom(q_{j})=\emptyset, i\neq j$. Given $d:\mathbb{R}\rightarrow
\mathbb{R}$, the evolution in time of $\mathcal{H}$ is described by the pair
of functions $x:\mathbb{R}\rightarrow X$, $q:\mathbb{R}\rightarrow Q$, where
$x(t)$ is the solution of the equation%
\begin{equation}
\dot{x}(t)=f_{\boldsymbol{I}(x(t))}(x(t),d(t))
\end{equation}%
with initial state $x(0)=x_{0}$ and%
\begin{equation}
q(t)=q_{\boldsymbol{I}(x(t))}
\end{equation} 
%
Let $p^{n}(t)$ and $v^{n}(t)$ denote the position on a horizontal axis and the velocity of vehicle $n$, respectively. The continuous state of $n+1$ is
\begin{equation}
{x}^{n+1}(t)=\left[
\begin{array}
[c]{c}%
p^{n}(t)-p^{n+1}(t)\\
v^{n}(t)-v^{n+1}(t)\\
v^{n}(t)
\end{array}
\right]
\end{equation}
\vspace{-0.1cm}
\begin{figure}[h!]
	\centering\includegraphics[width=1\columnwidth]{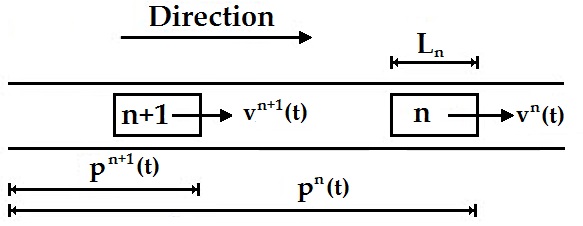}
\vspace{-0.7cm}
	\caption{The reference framework.}\label{Figure_reference}
\end{figure}%
\vspace{-0.1cm}
For simplicity, the dependency on $n$ will be omitted in the notation ${x}(t)
$. Moreover%
\begin{equation}
f_{i}(x(t),d(t))=\left[
\begin{array}
[c]{ccc}%
0 & 1 & 0\\
0 & 0 & 0\\
0 & 0 & 0
\end{array}
\right]  {x}(t)+\left[
\begin{array}
[c]{c}%
0\\
-1\\
0
\end{array}
\right]  g_{i}(x(t))+\left[
\begin{array}
[c]{c}%
0\\
1\\
1
\end{array}
\right]  d(t)
\end{equation}
where $u(t)=g_{i}(x(t))$, with $\left\vert u(t)\right\vert \leq a_{\max}$, is
the bounded state feedback control input and $d(t)$ is the control input of
the vehicle ahead. Since it is not known at $n+1$, it is modeled as a bounded
disturbance, i.e. $\left\vert d(t)\right\vert \leq a_{\max}$. The velocities are bounded too, i.e. $0\leq v^{n}(t)\leq v_{\max}$, $\forall \:
n=1,...,N$.

We define \emph{collision} the event that the distance between two vehicles is
less than $s_{n}$, which is the sum of the ahead vehicle length $L$ and a
minimum distance $L_{0}$. The functions $T_{E}:X\rightarrow\mathbb{R}$,
$T_{E}(x)=\frac{\left\vert x_{2}\right\vert }{a_{\max}}$, $T_{R}%
:X\rightarrow\mathbb{R}$, $T_{R}(x)=\frac{\left\vert x_{3}-x_{2}\right\vert
}{a_{\max}}$, and $T_{S}:X\rightarrow\mathbb{R}$, $T_{S}(x)=\lambda
\frac{\left\vert x_{3}-x_{2}\right\vert }{a_{\max}}$, represent, respectively,
the time headways needed to stop the vehicle starting from initial speed $x_{2}$ and
$\left\vert x_{3}-x_{2}\right\vert $, with deceleration $u(t)=-a_{\max}$, and
the time needed to stop the vehicle starting from initial speed $\left\vert
x_{3}-x_{2}\right\vert $, with $u(t)=-\frac{a_{\max}}{\lambda}$, $\lambda>1$.
We define the following thresholds for the space headway $p^{n}(t)-p^{n+1}(t)$:
\begin{itemize}
\item \emph{emergency distance} $\Delta E:\mathbb{R}^{3}\rightarrow
\mathbb{R}$
\footnotesize
\begin{equation}
\Delta E(x)=\left\{
\begin{array}
[c]{l}%
s_{n}\\
s_{n}+\frac{1}{2}a_{\max}T_{E}^{2}(x)
\end{array}
\right.
\begin{array}
[c]{c}%
x_{2}>0\\
x_{2}\leq0
\end{array}
\label{Delta_E}%
\end{equation}
\normalsize
represents the minimum distance where safety is ensured (see \cite{Gipps1981}). If at time $t$ the headway is equal to $\Delta E(x(t))$ and the leader starts braking with the maximum deceleration, provided that the follower also starts braking at the same time with the maximum deceleration, then collision is avoided.
\item \emph{risky distance} $\Delta R:\mathbb{R}^{3}\rightarrow\mathbb{R}$%
\footnotesize
\begin{equation}
\Delta R\left(  x\right)  =\left\{
\begin{array}
[c]{l}%
s_{n}+c_{r}T_{R}(x)x_{3}\\
s_{n}+c_{r}T_{R}(x)x_{3}+\frac{1}{2}a_{\max}T_{E}^{2}(x)
\end{array}
\right.
\begin{array}
[c]{c}%
x_{2}>0\\
x_{2}\leq0
\end{array}
\label{Delta_R}%
\end{equation}
\normalsize
has the same interpretation of the distance $\Delta E\left(  x\right)  $, but
it takes into account a human time-response, modeled by the add-on value
$c_{r}T_{R}(x)x_{3}$, where $c_{r}>0$ is a constant multiplication factor.
 Depending on the environment information
and on the human perception, such value can increase (more caution behaviour),
or decrease (more aggressive behaviour), but the condition $\Delta R\left(
x\right)  >\Delta E(x)$ is always satisfied.
\item \emph{safe distance} $\Delta S:\mathbb{R}^{3}\rightarrow\mathbb{R}$%
\footnotesize
\begin{equation}
\Delta S(x)=\left\{
\begin{array}
[c]{l}%
s_{n}+c_{s}T_{S}(x)x_{3}\\
s_{n}+c_{s}T_{S}(x)x_{3}+\frac{1}{2}a_{\max}T_{E}^{2}(x)
\end{array}
\right.
\begin{array}
[c]{c}%
x_{2}>0\\
x_{2}\leq0
\end{array}
\label{Delta_s}%
\end{equation}
\normalsize
considers a further safety margin w.r. to $\Delta R\left(  x\right)  $; in fact
$c_{s}\geq c_{r}$ 
 and $\lambda>1$ imply that $T_{S}>T_{R}$ and hence $\Delta S\left(x\right)>\Delta R\left(  x\right) $.

\item \emph{interaction distance}
$\Delta D:\mathbb{R}^{3}\rightarrow\mathbb{R}$%
\footnotesize
\begin{equation}
\Delta D\left(  x\right)  =\left\{
\begin{array}
[c]{l}%
s_{n}+c_{s}T_{S}(x)x_{3}\\
s_{n}+T_{D}(x_{3}-x_{2})
\end{array}
\right.
\begin{array}
[c]{c}%
x_{2}>0\\
x_{2}\leq0
\end{array}
\label{Delta_D}%
\end{equation}
\normalsize
where $T_{D}$ is a fixed time: it is the time headway beyond which one can be consider itself as a leader (see \cite{Wiedemann1991}, \cite{Fritzsche1994}).

Notice that when $x_{2}>0$, i.e. when $v^{n}(t)>v^{n+1}(t)$, $\Delta
D(x)=\Delta S(x)$.
\item \emph{approaching distance} $\Delta C:\mathbb{R}^{3}\rightarrow
\mathbb{R}$%
\footnotesize
\begin{equation}
\Delta C(x)=\left\{
\begin{array}
[c]{l}%
s_{n}+c_{s}T_{S}(x)x_{3}\\
s_{n}+ c_{s}T_{S}(x)x_{3} + c_{c} \mid\sqrt{-x_2}\mid
\end{array}
\right.
\begin{array}
[c]{c}%
x_{2}>0\\
x_{2}\leq0
\end{array}
\label{Delta_C}%
\end{equation}
\normalsize
is a threshold where the driver is approaching to high speed differences at short, decreasing distances (see \cite{Wiedemann1991}).
\end{itemize}
\begin{figure}[h!]
	\centering\includegraphics[width=1.1\columnwidth]{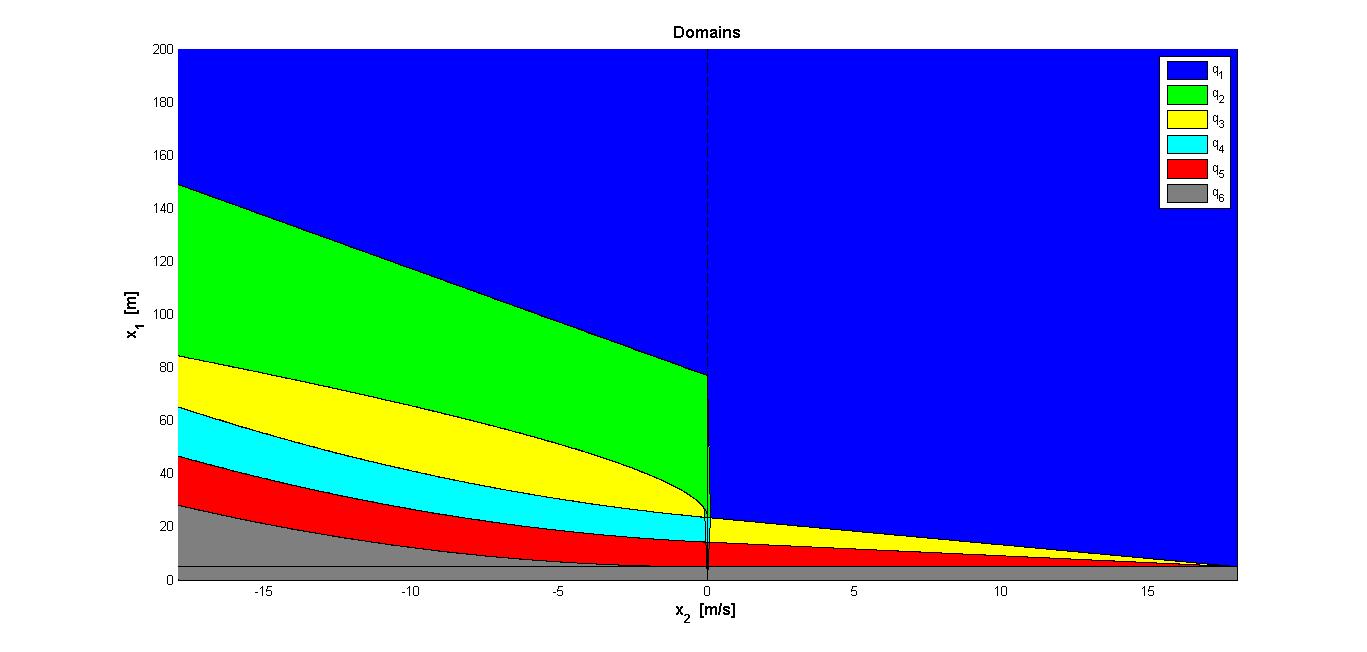}
\vspace{-0.8cm}
	\caption{The different thresholds and the domains defined by them at a fixed $\mathbf{v}^{n}=18$ m/s: $q_1$ domain is blue, while $q_2$ green, $q_3$ yellow, $q_4$ cyan, $q_5$ red and the unsafe zone $q_6$ is grey.}\label{Figure_DeltaX_V_diagram}
\end{figure}%
%
Setting $v^{n}(t)$ equal to a constant $\mathbf{v}^{n}$ (i.e. setting $d(t)=0$), in Figure \ref{Figure_DeltaX_V_diagram} the thresholds are represented in the bidimensional space, where on the horizontal axis is represented the speed difference $\mathbf{v}^{n}-v^{n+1}(t)=x_{2}(t)$, and on the vertical one the distance $p^{n}(t)-p^{n+1}(t)=x_{1}(t)$.
\begin{figure}[]
	\centering\includegraphics[width=0.9\columnwidth]{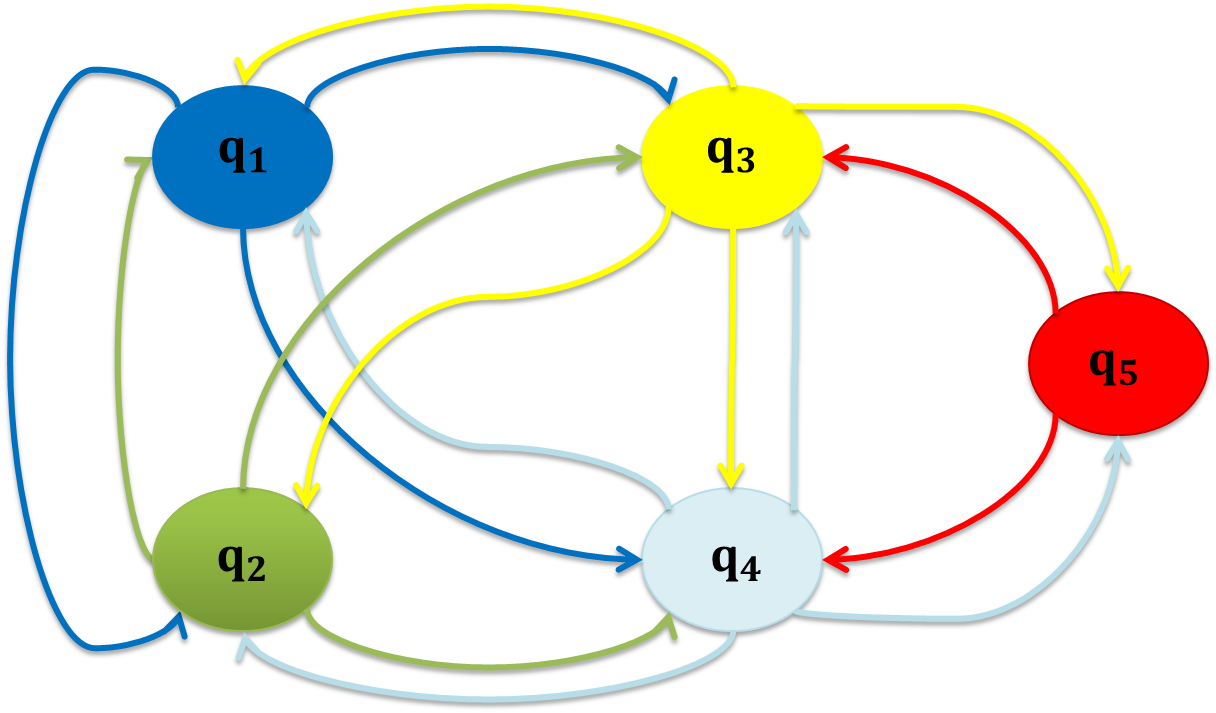}
\vspace{-0.2cm}
	\caption{The discrete states considered and their transitions. Colors are related to the domains depicted in Fig. \ref{Figure_DeltaX_V_diagram}.}\label{Figure_discrete_transitions}
\end{figure}%
%

The introduced thresholds define the domains of the continuous state space to be
associated to the discrete states. In the following description $\alpha_{1}$,
$\alpha_{2}$, $\alpha_{4}$ are positive sensitivity parameters, while
$v_{des}$ is the desired speed the driver wants to achieve. We assume that
$v_{des}=\mathbf{v}_{\max}$.%
\begin{enumerate}%
\item {$q_{1}$: Free driving}. $Dom(q_{1})$ is the set
\begin{align}
&\left\{  x\in X:(x_{1}>\Delta S(x))\wedge(x_{2}\geq0)\right\} \\
&\nonumber \cup\left\{x\in X:(x_{1}>\mathbf{m})\wedge(x_{2}<0)\right\} \label{Dom_q1}%
\end{align}
where $\mathbf{m}=\max\left\{  \Delta D(x),\Delta S(x)\right\}  $. The vehicle
can run freely because the leader vehicle is either too far away or faster or
both.
\begin{equation}
u_{n+1}(t)=g_{1}(x(t))=\alpha_{1}(v_{des}-\left(  x_{3}(t)-x_{2}(t)\right)
\label{u_q1}%
\end{equation}\label{q1control_action}%
%
\item {$q_{2}$: Following I}. $Dom(q_{2})$ is the set%
\begin{equation}
\left\{ x\in X:(x_{2}<0)\wedge\\
\nonumber (\mathbf{n} <x_{1}\leq\Delta D(x))\right\}
\end{equation}\label{Dom_q2}
with $\mathbf{n}=\max\left\{\Delta S(x),\Delta C(x)\right\}$.
Here the driver is closing in on the vehicle. The input depends on relative
speed and distance, following a modified version of the model in
\cite{Gazis1961}, i.e. 
\footnotesize
\begin{equation}
u_{n+1}(t)=g_{2}(x(t))=\alpha_{2}\cdot\frac{\left(  v_{des}+x_{2}(t)\right)
}{G-x_{1}(t)}\left(  x_{3}(t)-x_{2}(t)\right) \label{u_q2}%
\end{equation}\label{q2control_action}%
\normalsize
Here $G$ is a distance, and $G-x_{1}(t)>0$, $\forall \: x \in Dom(q_{2})$.%
\item {$q_{3}$: Following II}. $Dom(q_{3})$ is the set%
\begin{align}
&\left\{x\in X:(x_{2}\leq0)\wedge(\Delta S(x)< x_{1}<\mathbf{p})\right\}\cup \\
&\nonumber \left\{  x\in X:(x_{2}>0)\wedge\left(  \Delta R(x)<x_{1}\leq\Delta
S(x)\right)  \right\}%
\end{align}\label{Dom_q3}
where $\mathbf{p=}\min\left\{\Delta D(x),\Delta C(x)\right\}$. The
follower does not take any action, either because its speed is close to the leader's
one and the distance is small or because the speed difference is too large with respect to the distance.
\begin{equation}
u_{n+1}(t)=g_{3}(x(t))=0\label{u_q3}%
\end{equation}\label{q3control_action}%
%
\item {$q_{4}$: Closing in}. $Dom(q_{4})$ is the set%
\begin{align}
&\left\{  x\in X:(x_{2}\leq0)\wedge(\Delta R(x)<x_{1}\leq\Delta S(x))\right\}\\
\nonumber &\cup \left\{  x\in X:(x_{2}=0)\wedge(x_{1}=\Delta R(x))\right\} \label{Dom_q4}%
\end{align}\label{Dom_q4}%
the speed difference is large and the distance is not, so the driver has to
decelerate; he will do it depending on distance and relative speed, according
to the model in \cite{Saldana2000}.
\footnotesize
\begin{align}
&u_{n+1}(t)=g_{4}(x(t))= \\
\nonumber &\min\left\{  -\alpha_{4}\cdot\frac{x_{3}^{2}
(t)-(x_{3}(t)-x_{2}(t))^{2}}{2\left(  x_{1}(t)+s_{n}+\frac{c_{s}\lambda}%
{a_{max}}x^{2}_{3}(t)\right)  },\varepsilon\ast sign(x_{2})\right\} \label{u_q4}%
\end{align}\label{q4control_action}
\normalsize
%
The positive parameter $\varepsilon$ is introduced in order to present finite time convergence to the equilibrium points
, $\varepsilon\in\left(0, \frac{a_{max}}{\lambda}\right)$.%
\item {$q_{5}$: Danger}. $Dom(q_{5})$ is the set%
\begin{equation}
\begin{aligned}
&\left\{  x\in X:\left(  \Delta E(x)\leq x_{1}\leq\Delta R(x)\right)  \right\}\\
\nonumber &\backslash \left\{  x\in X:(x_{2}=0)\wedge(x_{1}=\Delta R(x))\right\}
\label{Dom_q5}%
\end{aligned}
\end{equation}
The distance from the $n$ vehicle is close to the unsafe one and the driver
uses his maximum deceleration.
\begin{equation}
u_{n+1}(t)=g_{5}(x(t))=-a_{\max}\label{u_q5}%
\end{equation}\label{q5control_action}%
\item {$q_{6}$: Unsafe}. Collision cannot be avoided: $Dom(q_{6})$ is the set%
\begin{equation}
\left\{  x\in X:x_{1}<\Delta E(x)\right\} \label{Dom_q6}%
\end{equation}%
\end{enumerate}
%
%
Finally, let us define the $Init$ set.%
\footnotesize
\begin{equation}
Init=\left(  \bigcup_{i=1}^{5}\{q_{i}\}\times\{Dom(q_{i})\cap\Sigma\}\right)
\label{Init_qxDom}%
\end{equation}
\normalsize
where%
\footnotesize
\begin{equation}%
\Sigma=\left\{  x:\left[
\begin{array}
[c]{c}%
s_{n}\\
-v_{\max}\\
0
\end{array}
\right]  \leq x\leq\left[
\begin{array}
[c]{c}%
\Delta_{\max}=500\\
v_{\max}\\
v_{\max}%
\end{array}
\right]  \right\} \label{Sigma}%
\end{equation}%
\normalsize
Notice that the unsafe domain is not included in the $Init$ set. Fig \ref{Figure_discrete_transitions} shows the admissible transitions.

For $n=1$, the first vehicle (the leader) can be described as an automaton
$\mathcal{H}_{1}$, obtained from $\mathcal{H}_{n+1}$ by setting $v^{0}%
(t)=v_{\max}$ and $Init=\{q_{1}\}\times\{Dom(q_{1})\cap\Sigma\}$. Therefore,
the hybrid state of the leader belongs to $Init$, until it remains leader,
i.e. until an ahead vehicle is not sufficiently close.
\begin{prop}
\label{p1}Let us assume $d:\mathbb{R}\rightarrow\mathbb{R}$ is given. The hybrid
automaton $\mathcal{H}_{n}$, $n=1,...,N$, is non-blocking, deterministic and non Zeno.
\end{prop}%
%
\begin{pf}
Non-blocking and determinism comes by construction. Moreover, by inspection,
any execution corresponding to
constant in time speed of the ahead vehicle is non Zeno. Moreover Zeno
executions are excluded if such speed is not constant, because of the bounds in
the acceleration.
\end{pf}%
Let us now analyze stability of $\mathcal{H}_{n}$. Suppose that $d(t)=0$,
and let $\mathbf{v}^{n}$ be given. Then the set of equilibrium states is
\begin{equation}
X_{e}=\left\{  x\in X:x_{2}=0 \wedge \Delta R(x)\leq x_{1}\leq\Delta S(x)\right\}
\end{equation}%
%
Let $x^{n}$ denote the state of $\mathcal{H}_{n}$, $n=1,...,N$ and $d^n(t)$ the associated disturbance.
Let us consider the compact set $\Omega=\Sigma\cap\bigcup_{i=1}^{5}Dom(q_{i})$, where $\Sigma$ has been defined
in (\ref{Sigma}). By construction such set is invariant for $\mathcal{H}_{n}$, i.e $\forall \: x^n_0 \in  \:\Omega, x^n(t) \in \Omega  \: \forall  \:t\geq0 \: \forall  \:d^n(t)$.%
\begin{prop}\label{stability}
Let us suppose that $d^n(t) = 0 \:\forall \:t\geq0$: then $\forall \: \varepsilon\in\left(0, \frac{a_{max}}{\lambda}\right)$ there exists $\widehat{t}$ such that for any initial state $x^n_{0}\in\Omega$, $x^n(\widehat{t})\in X_{e}$.
\end{prop}%
%
%
Let $\varepsilon$ and $\widehat{t}$ 
be given. Then we can analyze the platoon.
\begin{prop}
$\forall  \: x_{0}^{n}\in\Omega$, $x^{n}(t)\in X_{e}$, $\forall  \: t\geq
n\widehat{t}$, $\forall  \: n=1,...,N$.
\end{prop}%
Therefore stability is assured for a finite $N$.
%
%
%
\section{Mesoscopic model}\label{Mesoscopic model}
In real life microscopic model parameters are related to macroscopic quantities, such as traffic density. Being variance a density-dependent function (see \cite{Helbing1999}), a variance-driven adaptation mechanism is  adopted for changing thresholds depending on the local mean speed $\bar{v}_{n}$ value and local variance $\theta_{n}$ to improve the overall system performance.
In \cite{Helbing2006}, the authors formulate a variance-driven time headways
(VDT) model in terms of a meta-model to be applied to any car-following model
where a time headway $T_{0}$ can be expressed by a model parameter or a
combination of model parameters. They define a multiplication factor $\alpha_{T}(t)\in\lbrack1+\gamma V_{n}(t),\alpha_{T}^{\max}]$ for the time headways: $\alpha_{T}^{\max}$, $\gamma$ utilized parameters can be determined
from empirical data of the time-headway distribution for free and congested
traffic (see \cite{Helbing2006}), while $V_{n}$ is the variation coefficient
and is a function of mean speed and variance. With this formulation, they
relate on the driver acting not only to his own leader, but also to the
neighboring environment.
In this section, a similar mechanism is introduced, in order to allow the same possibilities. Furthermore, the case $\alpha_{T}<1$ will be considered as well, while maintaining the safety property. We make the assumption that $n+1$ vehicle receives information from the ahead vehicles
about their states.
\subsection{Variance-driven time headways model}\label{Variance-driven time headways model}
Inspired by \cite{Helbing2006}, in defining $\mathcal{H}_{n+1}$ a VDT-like factor $\alpha_{T}$ is obtained as
\begin{equation}
\alpha_{T}(t)=1+sat(z(t))\label{alfaT}  \: \: \: \: \: \: \alpha_{T}(t)\in\lbrack0,\alpha_{T}^{\max}]
\end{equation}%
where $z(t)$ is equal to %
\begin{equation}
\int_{0}^{t}\left[  -z(\tau)+\gamma V_{n}(\tau) sign\left(  \left(
x_{3}(\tau)-x_{2}(\tau)\right)  -\bar{v}_{n}(\tau)\right)  \right]
d\tau
\end{equation}\label{zetaFORalfa}%
with $V_{n}(t)$ being the variation coefficient,
\footnotesize
\begin{equation}
V_{n}(t) = \frac{\sqrt[2]{\theta_{n}(t)}}{\bar{v}_{n}(t)} = \frac{\sqrt[2]{\frac{1}{n}\sum_{i=1}^{n}\left(  v^{i}(t)-\left(  \frac{1}{n}\sum_{i=1}^{n}v^{i}(t)\right)
\right)^{2}}}{\frac{1}{n}\sum_{i=1}^{n}v^{i}(t)}
\end{equation}
\normalsize
%
and $\gamma$ the sensitivity of the time headway for velocity variations. The parameter defined in (\ref{alfaT}) allows to take into account both increasing and decreasing velocity variations.
Bounds on $\alpha_{T}(t)$ are used to define bounds on $z(t)$ with the saturation function $sat:\mathbb{R}\rightarrow\mathbb{R}$:
\begin{equation}
sat(z(t))=\left\{
\begin{array}
[c]{l}%
z_{T}^{\min}\\
z(t)\\
z_{T}^{\max}
\end{array}
\right.
\begin{array}
[c]{c}%
 z(t)\leq z_{T}^{\min}\\
\: z_{T}^{\min}<z(t)<z_{T}^{\max}\\
z(t)\geq z_{T}^{\max}
\end{array}
\label{Delta_E}%
\end{equation}
The positive saturation value will be $z_{T}^{\max}=\alpha_{T}^{\max}-1$, while the negative one will be $z_{T}^{\min}=-1+\alpha^{0}_{T}$, where $\alpha^{0}_{T}$ is a constant positive value for taking into account computation time.
%
%
%

The $\Delta R(x)$, $\Delta S(x)$, $\Delta C(x)$, $\Delta D(x)$ thresholds used for defining $\mathcal{H}_{n+1}$ will be modified by the $\alpha_T$ parameter altering $T_R(x)$, $T_S(x)$  and $T_D$. The resulting  $\Delta R'(x)$, $\Delta S'(x)$, $\Delta C'(x)$, $\Delta D'(x)$ thresholds will be defined by $T_R'(x)=\alpha_T T_R(x)$, $T_S'(x)=\alpha_T T_S(x)$ functions and $T_D'=\alpha_T T_D$. $\Delta E(x)$ will not be affected by the multiplication factor because of what it represents. It needs to be noticed that when $\alpha_{T}(t)=0$ the $\Delta R(x)$, $\Delta S(x)$ functions collapse in $\Delta E(x)$ one. Then, according to domain definitions, in the worst case considered safety is still ensured. It is possible to prove that Proposition \ref{stability} still holds for the modified model.
The $\mathcal{H}_{n+1}$ model is now able to take into account the neighboring environment: such possibility will suggest an anticipatory action to the driver (if ADAS case) or provide it autonomously (if ACC case).

In order to calculate the $\alpha_T$ parameter, on-line calculation of $\bar{v}_{n}$ and $\theta_{n}$ is expected: the needed information will be propagated by a vehicular network. Each vehicle will send data through such network: data will regard vehicle position and speed. Those data have to be received in "real-time" (compared to the human response time and safety-critical time-response) in order to select the right control action. In Section \ref{vehicular_network} a feasibility study of the needed communication network is described, which will allow us to consider such hypothesis as feasible.


%

%
\section{Vehicular networks}\label{vehicular_network}
Nowadays connected-vehicles (see \cite{Uhlemann}), namely vehicles that include
interactive advanced driver-assistance systems (ADASs) and cooperative
intelligent transport systems (C-ITSs), are a reality and in the next years
will be part of our daily life. 
%
Much interest is arising on C-ITSs, where vehicles cooperate by exchanging
messages wirelessly to achieve a higher level of safety and to avoid
on-the-road hazards. The cellular network 
is not a suitable choice
for safety applications due to stringent requirements for both bounded delay
and high reliability; consequently dedicated protocol stacks have been developed.
In this section we give a brief look at vehicular technologies and standards that can allow the operation of the proposed solution.
\subsection{Vehicular technologies}\label{vehicular_tecno}
WAVE (Wireless Access in Vehicular Environment) is the protocol stack defined
by the IEEE with the intent of extending the 802.11 family to include
vehicular environments. The early standards were approved in 2010 and,
commonly, WAVE refers to the IEEE 802.11p and IEEE 1609.x standards.

For the purpose of our application, it is important to remark the following
characteristics of IEEE 802.11p:
\begin{itemize}
  \item operation range up to 1000 m;
  \item communications in high-speed and high-mobility scenarios;
  \item priority and power control.
\end{itemize}
Among all supported architectures, the Independent BSS (IBSS) network topology
allows a set of stations to directly communicate with each other. This
capability, also known as ad hoc networking, achieves connectivity everywhere
since it does not require a fixed infrastructure to establish the connections.

Safety applications do not require all the features defined in the TCP/IP
network and transport layers that would also introduce unwanted overhead and
delay; to support them, the WAVE Short-Message Protocol (WSMP) was defined and
standardized through IEEE 1609.3. WSMP messages are allowed in the WAVE
Control Channel and therefore can benefit from a higher transmission power and
a favorable scheduling time.
\subsection{The network delay problem}\label{network_delay_problem}
Delay is a crucial parameter in real-time networks and services. The ITU-T
G.114 recommendation states that with a unidirectional end-to-end delay lower
than 150 ms, most applications, both speech and non-speech, will experience
essentially transparent interactivity.
For road safety applications, it is conventionally assumed that the maximum
end-to-end delay must be below 100 ms, while for traffic efficiency the
threshold rises to 500 ms.
On the other hand, vehicle drivers have different reaction times, depending on
specific stimuli. In \cite{Triggs1982} most unalerted drivers have shown
themselves capable of responding in less than 2.5 s in urgent situations.
The reaction time is lower if the stimuli is not only visual, but it is also vibratory or auditory
(see \cite{Pratichizzo2013}), so an ADAS can inform the driver of an imminent
hazard, to increase the probability of an in-time response.

Simulations conducted in \cite{Riihiijarvi} confirm that the Control Channel
has good characteristics in terms of End-to-End Delay and Radio Range for
safety messages (lower than 20 ms and higher than 750 m, respectively) when
the message rate is lower than 1000 packets/s and the message size is below
256 Bytes. Estimating the delay in a multi-hop scenario is not trivial, since
many variables (e.g., distance, interference, computing time) impact on it,
but with some assumptions we can respect an upper-bound value of 100 ms. The
proposed controller considers a vehicle as a predecessor or a follower only if
its distance is up to 500 m but the WAVE technology allows to send data more
than 750 m away (in \cite{Riihiijarvi} the maximum distance reached is around
2.5 km in open air scenarios). Consequently the number of relays needed to
reach the last vehicle is not strictly bounded to the number of vehicles in
the queue but depends essentially on the distance. Considering clusters up to
5 vehicles, the worst case occurs when vehicles are along a line 500 m apart
from each other and the channel quality allows a maximum range of 750 m: to
reach the end of the line under such conditions, 4 relays are needed with a
total delay lower than 80 ms plus the eventual computation time on the relay
nodes. As the environment is highly dynamic, the topology changes over time
but, if the vehicles are still clustered, the delay and the number of hops are
surely better than in the worst case. It follows that the reduction of the
number of vehicles in a single cluster and the minimization of the processing
operations on the information to be relayed are sufficient to comply the
constraints of safety applications.

It is worth noting that, when more vehicles are within the Radio Range, the
same message may arrive several times to the same vehicle. Moreover, the
message arrival does not implicitly provide information about the positions of
vehicles that either generated or relayed it. The first problem is easily
solved including Timestamp and Serial Number in the packet in order to enable
the Duplicate Discovery Process; the second issue can be simply managed using
GPS geolocation.
%
%
\section{Simulation results} \label{Simulation results}
In this section simulation results for the developed hybrid model are provided in Simulink. 
The targets we wanted to achieve are to verify the ability to correctly represent a safe car-following situation (collision avoidance) and to compare the use or disuse of the VDT-like mechanism in a complex traffic scenario.
%
The scenario proposed in order to validate results is composed by 5 vehicles, whose parameters are described in \cite{ADHS_Iovine_Arxiv}. Each vehicle is supposed to be driven by an ACC system controlled by the introduced mesoscopic hybrid automaton.

At the beginning, the first vehicle is the leader of a platoon composed by 4 vehicles: at time 30" it decelerates for reaching the speed of 18 m/s, while at time 90" it speeds up to a velocity of 33 m/s. The fifth vehicle reaches and tag along to the platoon: its initial separation is the maximum we consider for an interaction between two vehicles (segment road of 500 meters), according to Section \ref{vehicular_network}: the technological constraints are then fulfilled. The 5th vehicle is approaching a platoon that is decreasing its initial speed: the $\alpha_T$ parameter utilization is expected to anticipate its deceleration. Both the active or inoperative VDT-like mechanism mode are implemented and compared.
%
%
%

\begin{figure}[h!]
	\centering\includegraphics[width=1\columnwidth]{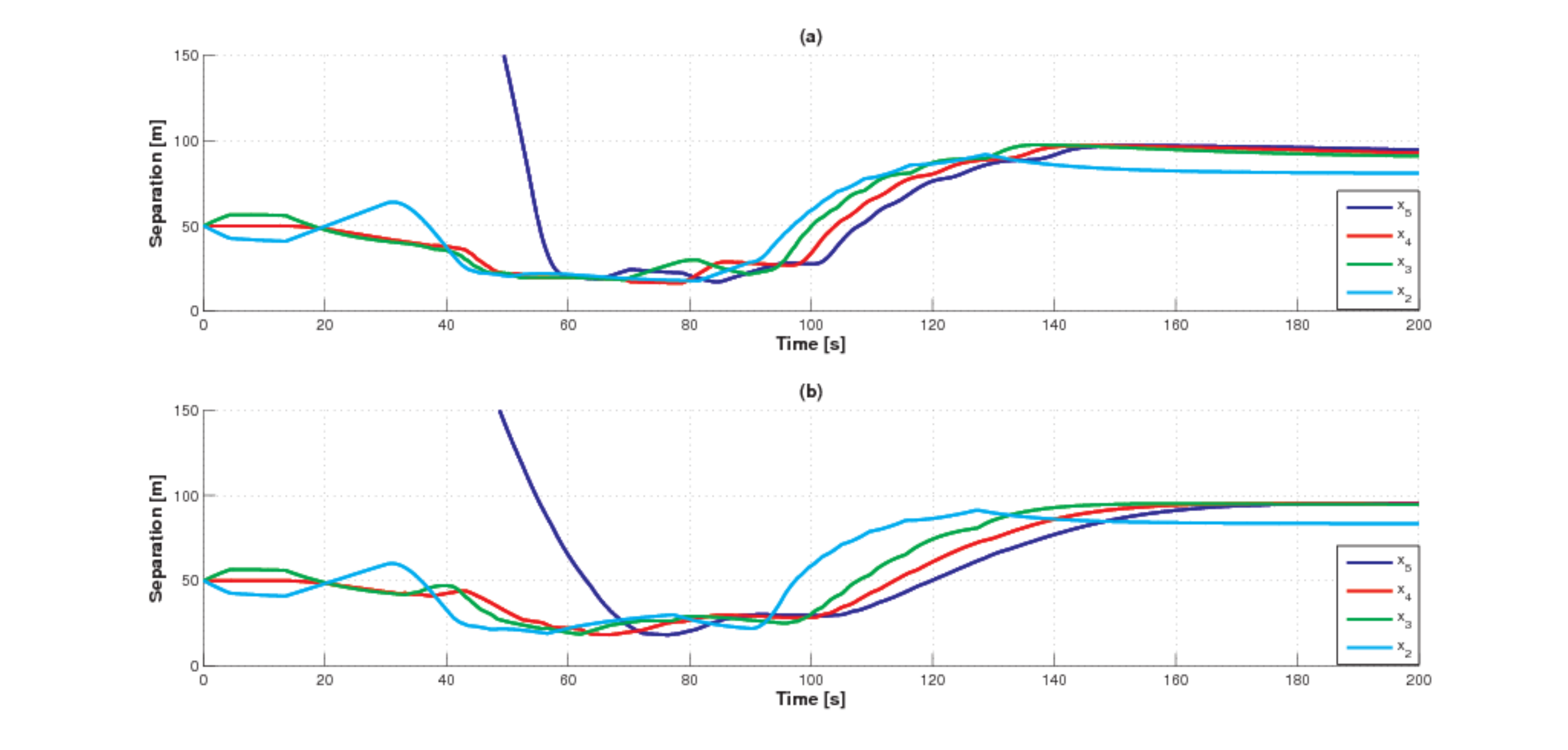}
\vspace{-0.7cm}
	\caption{The separations among vehicles in the active VDT-like mechanism case, (a), and the inoperative one, (b). The fifth vehicle dynamics is depicted in blue color.}\label{Figure_separation_simulation}
\end{figure}

\begin{figure}
	\centering\includegraphics[width=1.1\columnwidth]{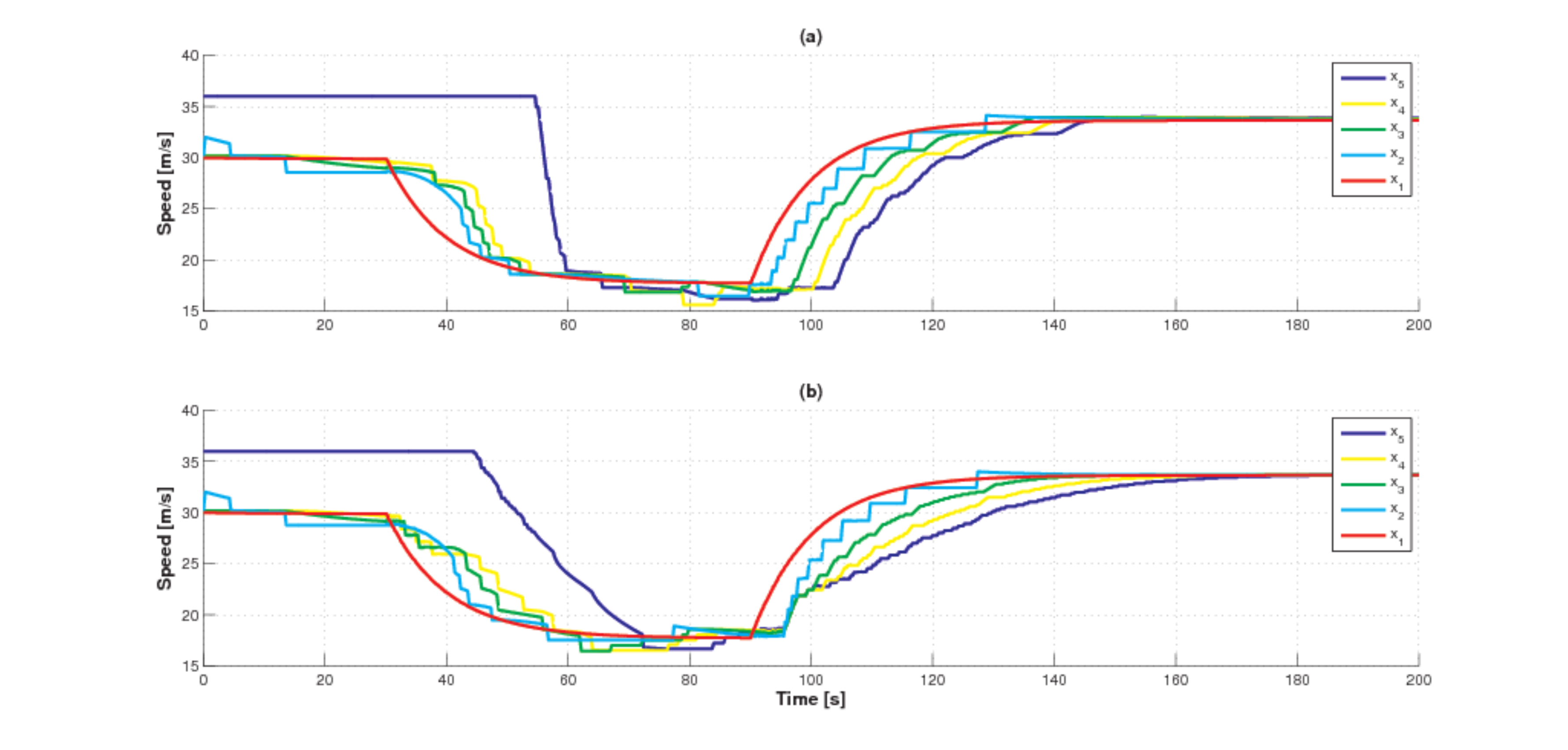}
\vspace{-0.7cm}
    \caption{The speed of each vehicle  in the active VDT-like mechanism case, (a), and the inoperative one, (b). The first vehicle dynamics is depicted in red color, while the fifth one is in blue.}\label{Fig_comparison_sino_alphat}
\end{figure}

\begin{figure}
	\centering\includegraphics[width=1.4\columnwidth]{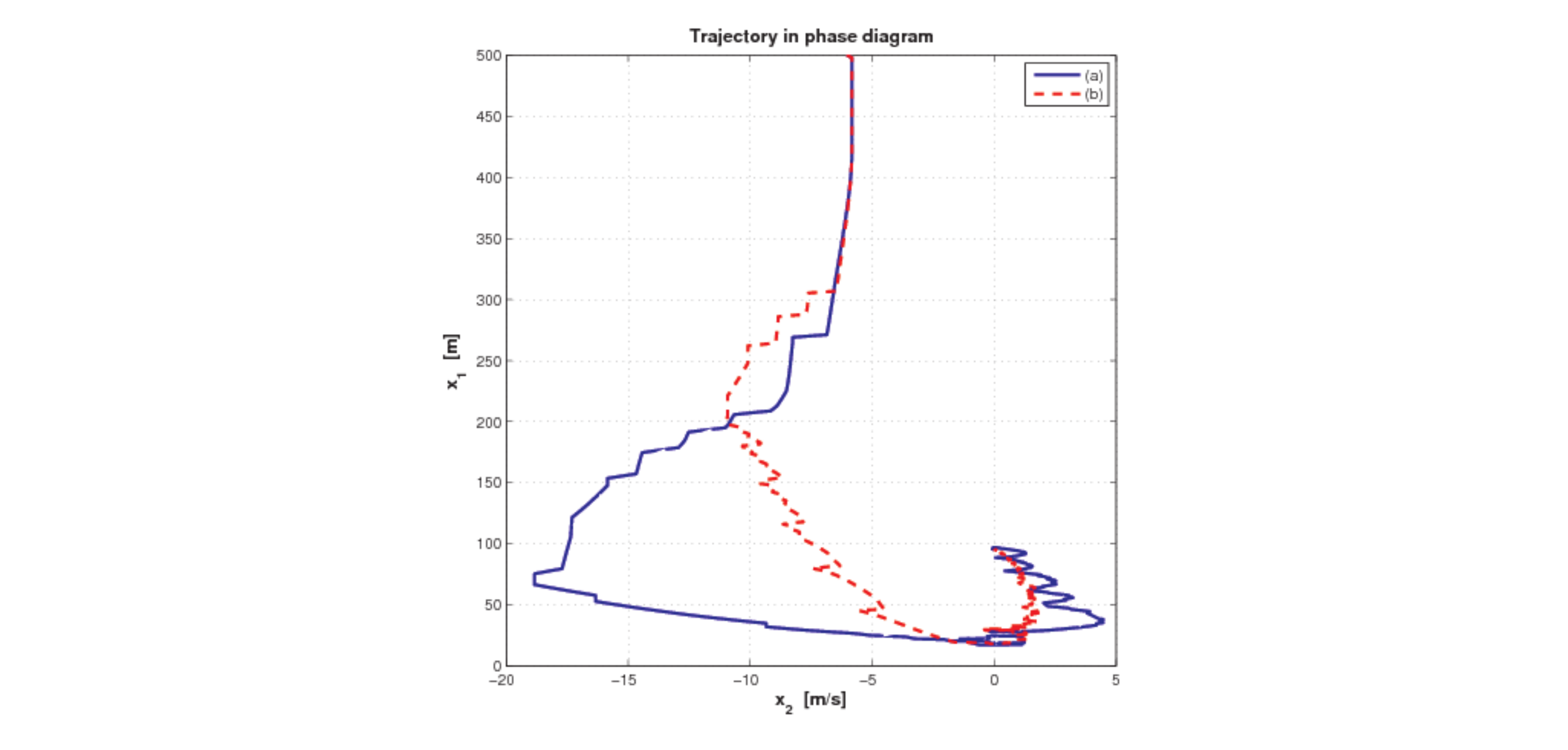}
\vspace{-0.7cm}
	\caption{Trajectories in phase portrait are depicted: in (a) (continue curve) the platoon does not use VDT-like mechanism, while in (b) (dotted curve) it does.}\label{Figure_phase_dynamics}
\end{figure}

Figure \ref{Figure_separation_simulation} shows the separations among vehicles: there is no event of collision for both cases.
In Figure \ref{Fig_comparison_sino_alphat} the anticipatory action of the VDT is clearly depicted: indeed, in (b) the fifth vehicle starts decelerating before than in (a), around 45" and 55", respectively, as can be stated by the blue curve.
A similar anticipatory action is seen when acceleration occurs: in (b), the fifth vehicle accelerates before the point representing 100", while in (a) it does after that point.
This will lead to a smoother behavior, which is described in the $\mathcal{H}_{5}$ model phase portrait in Figure \ref{Figure_phase_dynamics}. %
%
It is possible to verify some oscillations: they are expected around the equilibrium point because of the follower inability to precisely predict leader state. In the simulations, less oscillations appear in case of VDT-like mechanism utilization. The improvements are due to the use of more information. 
Those results are important because non-smooth transients are responsible for string instability and shock waves propagation. Hence, less oscillations and less magnitude of oscillations, as shown in our results, provide a better response for these problems. 
In Figure \ref{Figure_phase_dynamics} steady-state behaviours are close to each other, while the transient is smoother in case of VDT-like mechanism utilization. The simulations then show how it is possible to improve ACC performances taking into account information regarding the neighbors vehicles.
\section{Conclusions}\label{Conclusions}
In this paper, an innovative mesoscopic hybrid model for car-following situations is introduced: its purpose is to safely control the single vehicle dynamics through an automatic controller that would replace the human control actions or to support human driver with an assistance system human-inspired.
The controller processes other vehicles information and takes decisions about braking or throttle actions. 
It does not contain human weaknesses, such as mistakes or distractions, but strengths such as adaptability to various conditions and comfort.

In the hybrid automaton, microscopic models representing human-drivers imitate a safe human behavior; indeed, systems properties and a stability definition with respect to some desired equilibrium points have been investigated and proved. To better represent human-style, a macroscopic value dependence is added to the final model: the consideration of macroscopic quantities is an important step-ahead.
A feasibility study of the needed communication network is analyzed and described and simulation results confirm the necessity to consider the two description levels for improving traffic throughput.
%
%
%
\scriptsize
\bibliography{mcnbib}
\normalsize
%
\end{document}